\documentclass[12pt,reqno]{amsart}
\usepackage{amsmath, amsthm, amscd, amsfonts, amssymb, graphicx, color}
\usepackage[bookmarksnumbered, colorlinks, plainpages]{hyperref}
\DeclareMathOperator{\Tr}{Tr}
\theoremstyle{definition}

\theoremstyle{remark}

\numberwithin{equation}{section}
\begin{document}
\begin{center}
{\textbf{Conformal Yamabe soliton and $*$-Yamabe soliton with torse forming potential vector field }}
\end{center}
\vskip 0.3cm
\begin{center}By\end{center}\vskip 0.3cm
\begin{center}
{Soumendu Roy \footnote{The first author is the corresponding author, supported by Swami Vivekananda Merit Cum Means Scholarship, Government of West Bengal, India.}, Santu Dey $^2$ and Arindam~~Bhattacharyya $^3$}
\end{center}
\vskip 0.3cm
\address[Soumendu Roy]{Department of Mathematics,Jadavpur University, Kolkata-700032, India}
\email{soumendu1103mtma@gmail.com}

\address[Santu Dey]{Department of Mathematics, Bidhan Chandra College, Asansol - 4, West Bengal-713304 , India}
\email{santu.mathju@gmail.com}

\address[Arindam Bhattacharyya]{Department of Mathematics,Jadavpur University, Kolkata-700032, India}
\email{bhattachar1968@yahoo.co.in}
\vskip 0.5cm
\begin{center}
\textbf{Abstract}\end{center}
The goal of this paper is to study conformal Yamabe soliton and $*$-Yamabe soliton, whose potential vector field is torse forming. Here, we have characterized conformal Yamabe soliton admitting potential vector field as torse forming  with respect to Riemannian connection, semi-symmetric metric connection and projective semi-symmetric connection on Riemannian manifold. We have also shown the nature of $*$-Yamabe soliton with torse forming vector field on Riemannian manifold admitting Riemannian connection. Lastly we have developed an example to corroborate some theorems regarding Riemannian connection on Riemannian manifold.
\vspace {0.3cm}
\section{\textbf{Introduction}}
The concept of Yamabe flow was first introduced by Hamilton \cite{hamil} to construct Yamabe metrics on compact Riemannian manifolds. On a Riemannian or pseudo-Riemannian manifold $M$, a time-dependent metric $g(\cdot, t)$ is said to evolve by the Yamabe flow if the metric $g$ satisfies the given equation,
\begin{equation}\label{1.1}
  \frac{\partial }{\partial t}g(t)=-rg(t),\hspace{0.5cm} g(0)=g_{0},
\end{equation}
where $r$ is the scalar curvature of the manifold $M$.\\\\
In 2-dimension the Yamabe flow is equivalent to the Ricci flow \cite{rsham} (defined by $\frac{\partial }{\partial t}g(t) = -2S(g(t))$, where $S$ denotes the Ricci tensor). But in dimension $> 2$ the Yamabe and Ricci flows do not agree, since the Yamabe flow preserves the conformal class of the metric but the Ricci flow does not in general.\\
A Yamabe soliton \cite{barbosa} correspond to self-similar solution of the Yamabe flow, is defined on a Riemannian or pseudo-Riemannian manifold $(M, g)$ as:
\begin{equation}\label{1.2}
  \frac{1}{2}\pounds_V g = (r-\lambda)g,
\end{equation}
where $\pounds_V g$ denotes the Lie derivative of the metric $g$ along the vector field $V$, $r$ is the scalar curvature and $\lambda$ is a constant. Moreover a Yamabe soliton is said to be expanding, steady, shrinking depending on $\lambda$ being positive, zero, negative respectively. If $\lambda$ is a smooth function then \eqref{1.2} is called almost Yamabe soliton \cite{barbosa}. \\
Many authors such as \cite{roy2}, \cite{roy3}, \cite{dong}, \cite{ghosh} have studied Yamabe solitons on some contact manifolds.\\
In 2015, N. Basu and A. Bhattacharyya \cite{nbab} established the notion of conformal Ricci soliton \cite{soumendu}, \cite{roy} as:
\begin{equation}\label{1.3}
\pounds_V g + 2S = \Big[2\lambda - \Big(p + \frac{2}{n}\Big)\Big]g,
\end{equation}
where $S$ is the Ricci tensor, $p$ is a scalar non-dynamical field(time dependent scalar field),  $\lambda$ is constant, $n$ is the dimension of the manifold. \\
Using \eqref{1.2} and \eqref{1.3}, we introduce the notion of conformal Yamabe soliton as:\\\\
\textbf{Definition 1.1:} A Riemannian or pseudo-Riemannian manifold $(M,g)$ of dimension $n$ is said to admit conformal Yamabe soliton if
\begin{equation}\label{1.4}
 \pounds_V g+\Big[2\lambda -2r - \Big(p + \frac{2}{n}\Big)\Big]g=0,
 \end{equation}
 where $\pounds_V g$ denotes the Lie derivative of the metric $g$ along the vector field $V$, $r$ is the scalar curvature and $\lambda$ is a constant, $p$ is a scalar non-dynamical field(time dependent scalar field), $n$ is the dimension of the manifold. The conformal Yamabe soliton is said to be expanding, steady, shrinking depending on $\lambda$ being positive, zero, negative respectively. If the vector field $V$ is of gradient type i.e $V = grad(f)$, for $f$ is a smooth function on $M$, then the equation \eqref{1.4} is called conformal gradient Yamabe soliton. \\\\
 The notion of $*$-Ricci tensor on almost Hermitian manifolds and $*$-Ricci tensor of real hypersurfaces in non-flat complex space were introduced by Tachibana \cite{tachi} and Hamada \cite{hama} respectively where the $*$-Ricci tensor is defined by:
\begin{equation}\label{1.5}
  S^*(X,Y)=\frac{1}{2}(\Tr \{\varphi\circ R(X,\varphi Y)\}),
\end{equation}
for all vector fields $X,Y$ on $M^n$, $\varphi$ is a (1,1)-tensor field and $\Tr$ denotes Trace.\\
If $S^*(X,Y)=\lambda g(X,Y)+\nu \eta(X)\eta(Y)$ for all vector fields $X,Y$ and $\lambda$, $\nu$ are smooth functions, then the manifold is called $*$-$\eta$-Einstein manifold.\\
Further if $\nu=0$ i.e $S^*(X,Y)=\lambda g(X,Y)$ for all vector fields $X,Y$ then the manifold becomes $*$-Einstein.\\
In 2014, Kaimakamis and Panagiotidou \cite{kaipan} introduced the notion of $*$-Ricci soliton which can be defined as:
\begin{equation}\label{1.6}
  \pounds_V g + 2S^* + 2\lambda g=0,
\end{equation}
for all vector fields $X,Y$ on $M^n$ and $\lambda$ being a constatnt.\\
Using \eqref{1.2} and \eqref{1.6}, we develop the notion of $*$-Yamabe soliton as:\\\\
\textbf{Definition 1.2:} A Riemannian or pseudo-Riemannian manifold $(M,g)$ of dimension $n$ is said to admit $*$-Yamabe soliton if
\begin{equation}\label{1.7}
  \frac{1}{2}\pounds_V g = (r^*-\lambda)g,
\end{equation}
where $\pounds_V g$ denotes the Lie derivative of the metric $g$ along the vector field $V$, $r^* = \Tr(S^*)$ is the $*$- scalar curvature and $\lambda$ is a constant. The $*$-Yamabe soliton is said to be expanding, steady, shrinking depending on $\lambda$ being positive, zero, negative respectively. If the vector field $V$ is of gradient type i.e $V = grad(f)$, for $f$ is a smooth function on $M$, then the equation \eqref{1.7} is called $*$-gradient Yamabe soliton.\\\\
The outline of the article goes as follows:\\
In section 2, after a brief introduction, we have discussed some needful results which will be used in the later sections. Section 3 deals with some applications of torse forming potential vector field on conformal Yamabe soliton. In this section we have contrived conformal Yamabe soliton own up to Riemannian connection, semi-symmetric metric connection and projective semi-symmetric connection with torse forming vector field to accessorize the nature of this soliton on Riemannian manifold and we have proved Theorem 3.1, Theorem 3.3 and Theorem 3.5 concerning those mentioned connections. Section 4 is devoted to utilize of torse forming potential vector field on $*$-Yamabe soliton  with respect to Riemannian connection and we have evolved a theorem to develop the essence of this soliton. In section 5, we have constructed an example to illustrate the existence  of the conformal Yamabe soliton on 3-dimensional Riemannian manifold.\\\\
\section{\textbf{Preliminaries}}
A nowhere vanishing vector field $\tau$ on a Riemannian or pseudo-Riemannian manifold $(M,g)$ is called torse-forming \cite{yano} if
\begin{equation}\label{2.1}
  \nabla_X \tau = \phi X+\alpha(X)\tau,
\end{equation}
where $\nabla$ is the Levi-Civita connection of $g$, $\phi$ is a smooth function and $\alpha$ is a 1-form.
Moreover The vector field $\tau$ is called concircular \cite{chen}, \cite{kyano} if the 1-form $\alpha$ vanishes identically in the equation \eqref{2.1}. The vector field $\tau$ is called concurrent \cite{Schouten}, \cite{kyano1} if in \eqref{2.1} the 1-form $\alpha$ vanishes identically and the function $\phi = 1$. The vector field $\tau$ is called recurrent if in \eqref{2.1} the function $\phi = 0$. Finally if in \eqref{2.1} $\phi = \alpha = 0$, then the vector field $\tau$ is called a parallel vector field.\\
In 2017, Chen \cite{chen1} introduced a new vector field called torqued vector field. If the vector field $\tau$ staisfies \eqref{2.1} with $\alpha(\tau) = 0$, then $\tau$ is called torqued vector field. Also in this case, $\phi$ is known as the torqued function and the 1-form $\alpha$ is the torqued form of $\tau$.\\
From \cite{Friedmann}, \cite{Hayden}, \cite{kyano2}, the relation between the semi-symmetric metric connection $\bar{\nabla}$ and the connection $\nabla$ of $M$ is given by:
\begin{equation}\label{2.2}
  \bar{\nabla}_X Y = \nabla_X Y + \pi(Y)X - g(X, Y)\rho,
\end{equation}
where $\pi(X) = g(X, \rho), \forall X \in \chi(M)$, the Lie algebra of vector fields of $M$.\\
Also the Riemannian curvature tensor $\bar{R}$, Ricci tensor $\bar{S}$ and the scalar curvature $\bar{r}$ of $M$ associated with the semi-symmetric metric connection $\bar{\nabla}$ are given by \cite{Friedmann}:
\begin{equation}\label{2.2new}
  \bar{R}(X,Y)Z=R(X,Y)Z-P(Y,Z)X+P(X,Z)Y-g(Y,Z)LX+g(X,Z)LY,
\end{equation}
\begin{equation}\label{2.3}
  \bar{S}(Y, Z) = S(Y, Z) - (n - 2)P(Y, Z) - ag(Y, Z),
\end{equation}
\begin{equation}\label{2.4}
  \bar{r} = r - 2(n - 1)a,
\end{equation}
where $P$ is a (0,2) tensor field given by: $$P(X, Y) =g(LX,Y)= (\nabla_X \pi)(Y) - \pi(X)\pi(Y) +\frac{1}{2}\pi(\rho)g(X, Y), \forall X,Y \in \chi(M)$$ and $a=\Tr(P)$.\\
From \cite{Zhao}, the relation between projective semi-symmetric connection $\tilde{\nabla}$ and the connection $\nabla$  is given by:
\begin{equation}\label{2.5}
  \tilde{\nabla}_X Y = \nabla_X Y + \psi(Y)X + \psi(X)Y + \mu(Y)X - \mu(X)Y,
\end{equation}
where the 1-forms $\psi$ and $\mu$ are given by:
$$\psi(X) = \frac{n-1}{2(n+1)}\pi(X),$$
 $$\mu(X) = \frac{1}{2}\pi(X)$$.\\
Also the Riemannian curvature tensor $\tilde{R}$, Ricci tensor $\tilde{S}$ and the scalar curvature $\tilde{r}$ of $M$ associated with the projective semi-symmetric connection $\tilde{\nabla}$ are given by \cite{Zhao}, \cite{shaikh}:
\begin{equation}\label{2.5new}
  \tilde{R}(X,Y)Z=R(X,Y)Z+\theta(X,Y)Z+\omega(X,Z)Y-\omega(Y,Z)X,
\end{equation}
\begin{equation}\label{2.6}
  \tilde{S}(Y,Z)=S(Y, Z) + \theta(Y, Z) - (n - 1)\omega(Y, Z),
\end{equation}
\begin{equation}\label{2.7}
  \tilde{r}= r + \Tr(\theta) - (n - 1)\Tr(\omega),
\end{equation}
where
$$\theta(X,Y)=\frac{1}{2} [(\nabla_Y \pi)(X) - (\nabla_X\pi)(Y)],$$
$$\omega(X,Y)= \frac{n - 1}{2(n + 1)} (\nabla_X\pi)(Y)+\frac{1}{2}(\nabla_Y \pi)(X)-\frac{n^2}{(n+1)^2}\pi(X)\pi(Y),$$
$\forall X,Y,Z \in \chi(M).$
\vspace {0.3cm}
\section{\textbf{Application of torse forming vector field on conformal Yamabe soliton}}
Let $(g,\tau,\lambda)$ be a conformal Yamabe soliton on $M$ with respect to the Riemannian connection $\nabla$. Then from \eqref{1.4} we have,
\begin{equation}\label{3.1}
  (\pounds_\tau g)(X,Y)+\Big[2\lambda -2r - \Big(p + \frac{2}{n}\Big)\Big]g(X,Y)=0.
\end{equation}
Now using\eqref{2.1}, we obtain,
\begin{eqnarray}\label{3.2}
  (\pounds_\tau g)(X,Y) &=& g(\nabla_X \tau,Y)+g(X,\nabla_Y \tau) \nonumber\\
                        &=& 2\phi g(X,Y)+\alpha(X)g(\tau,Y)+\alpha(Y)g(\tau,X),
\end{eqnarray}
for all $X,Y \in M$.\\
Then using \eqref{3.2}, \eqref{3.1} becomes,
\begin{equation}\label{3.3}
  \Big[r-\phi-\lambda+\frac{1}{2}\Big(p + \frac{2}{n}\Big)\Big]g(X,Y)=\frac{1}{2}\Big[\alpha(X)g(\tau,Y)+\alpha(Y)g(\tau,X)\Big].
\end{equation}
Taking contraction of \eqref{3.3} over X and Y, we have,
\begin{equation}\label{3.4}
  \Big[r-\phi-\lambda+\frac{1}{2}\Big(p + \frac{2}{n}\Big)\Big]n=\alpha(\tau),
\end{equation}
which leads to,
\begin{equation}\label{3.5}
  \lambda=r-\phi-\frac{\alpha(\tau)}{n}+\frac{1}{2}\Big(p + \frac{2}{n}\Big).
\end{equation}
So we can state the following theorem:\\\\
\textbf{Theorem 3.1.} {\em Let $(g,\tau,\lambda)$ be a conformal Yamabe soliton on $M$ with respect to the Riemannian connection $\nabla$. Then the vector field $\tau$ is torse-forming if $\lambda=r-\phi-\frac{\alpha(\tau)}{n}+\frac{1}{2}(p + \frac{2}{n})$, is constant and the soliton is expanding, steady, shrinking according as $r-\phi-\frac{\alpha(\tau)}{n}+\frac{1}{2}(p + \frac{2}{n}) \gtreqqless 0$. }\\\\
Now in \eqref{3.5}, if the 1-form $\alpha$ vanishes identically then $\lambda=r-\phi+\frac{1}{2}(p + \frac{2}{n}).$\\\\
If the 1-form $\alpha$ vanishes identically and the function $\phi = 1$ in \eqref{3.5}, then $\lambda=r-1+\frac{1}{2}(p + \frac{2}{n}).$\\\\
In \eqref{3.5}, if the function $\phi = 0$, then $\lambda=r-\frac{\alpha(\tau)}{n}+\frac{1}{2}(p + \frac{2}{n}).$\\\\
If $\phi = \alpha = 0$ in \eqref{3.5}, then $\lambda=r+\frac{1}{2}(p + \frac{2}{n}).$\\\\
Finally in \eqref{3.5}, if $\alpha(\tau) = 0$, then $\lambda=r-\phi+\frac{1}{2}(p + \frac{2}{n}).$\\
Then we have,\\\\
\textbf{Corollary 3.2.} {\em Let $(g,\tau,\lambda)$ be a conformal Yamabe soliton on $M$ with respect to the Riemannian connection $\nabla$. Then the vector field $\tau$ is\\\\
(i)concircular if $\lambda=r-\phi+\frac{1}{2}(p + \frac{2}{n})$, is constant and the soliton is expanding, steady, shrinking according as $r-\phi+\frac{1}{2}(p + \frac{2}{n})\gtreqqless 0$.\\\\
(ii)concurrent if $\lambda=r-1+\frac{1}{2}(p + \frac{2}{n})$, is constant and the soliton is expanding, steady, shrinking according as $r-1+\frac{1}{2}(p + \frac{2}{n})\gtreqqless 0$.\\\\
(iii)recurrent if $\lambda=r-\frac{\alpha(\tau)}{n}+\frac{1}{2}(p + \frac{2}{n})$, is constant and the soliton is expanding, steady, shrinking according as $r-\frac{\alpha(\tau)}{n}+\frac{1}{2}(p + \frac{2}{n}) \gtreqqless 0$.\\\\
(iv)parallel if $\lambda=r+\frac{1}{2}(p + \frac{2}{n})$, is constant and the soliton is expanding, steady, shrinking according as $r+\frac{1}{2}(p + \frac{2}{n}) \gtreqqless 0$.\\\\
(v)torqued if $\lambda=r-\phi+\frac{1}{2}(p + \frac{2}{n})$, is constant and the soliton is expanding, steady, shrinking according as  $r-\phi+\frac{1}{2}(p + \frac{2}{n}) \gtreqqless 0$.}\\\\
Let us now consider $(g,\tau,\lambda)$ as a conformal Yamabe soliton on $M$ with respect to the semi-symmetric metric connection $\bar{\nabla}$. Then we have,
\begin{equation}\label{3.6}
  (\bar{\pounds}_\tau g)(X,Y)+[2\lambda -2\bar{r} - (p + \frac{2}{n})]g(X,Y)=0,
\end{equation}
where $\bar{\pounds}_\tau$ is the Lie derivative along $\tau$ with respect to $\bar{\nabla}$.\\
Now using \eqref{2.2}, we get,
\begin{eqnarray}\label{3.7}
 (\bar{\pounds}_\tau g)(X,Y) &=& g(\bar{\nabla}_X \tau,Y)+g(X,\bar{\nabla}_Y \tau) \nonumber\\
                             &=& g(\nabla_X \tau + \pi(\tau)X - g(X, \tau)\rho, Y)\nonumber\\
                             &+& g(X,\nabla_Y \tau + \pi(\tau)Y - g(Y, \tau)\rho)\nonumber\\
                             &=&(\pounds_\tau g)(X,Y)+2\pi(\tau)g(X,Y)\nonumber\\
                             &-&[g(x,\tau)\pi(Y)+g(Y,\tau)\pi(X)].
\end{eqnarray}
Using \eqref{3.2} in \eqref{3.7}, we obtain,
\begin{eqnarray}\label{3.8}
  (\bar{\pounds}_\tau g)(X,Y)&=& 2\phi g(X,Y)+\alpha(X)g(\tau,Y)+\alpha(Y)g(\tau,X) +2\pi(\tau)g(X,Y)\nonumber\\
                             &-& [g(x,\tau)\pi(Y)+g(Y,\tau)\pi(X)].
\end{eqnarray}
From \eqref{2.4} and \eqref{3.8}, \eqref{3.6} becomes,
\begin{multline}\label{3.9}
  \Big[\phi+\pi(\tau)-r+2(n-1)a+\lambda-\frac{1}{2} \Big(p + \frac{2}{n}\Big)\Big] g(X,Y)\\
  + \frac{1}{2}\Big[\Big\{\alpha(X)-\pi(X)\Big\} g(\tau,Y)+\Big\{\alpha(Y)-\pi(Y)\Big\} g(\tau,X)\Big] =0.
  \end{multline}
Taking contraction of \eqref{3.9} over X and Y, we have,
\begin{equation}\label{3.10}
  \Big[\phi-r+2(n-1)a+\lambda-\frac{1}{2} \Big(p + \frac{2}{n}\Big)\Big]n+(n-1)\pi(\tau)+\alpha(\tau)=0,
\end{equation}
which leads to,
\begin{equation}\label{3.11}
  \lambda=r-\phi-2(n-1)a+\frac{1}{2}\Big(p + \frac{2}{n}\Big)-\frac{n-1}{n}\pi(\tau)-\frac{\alpha(\tau)}{n}.
\end{equation}
Hence we can state the following:\\\\
\textbf{Theorem 3.3} {\em Let $(g,\tau,\lambda)$ be a conformal Yamabe soliton on $M$ with respect to the semi-symmetric metric connection $\bar{\nabla}$. Then the vector field $\tau$ is torse-forming if $\lambda=r-\phi-2(n-1)a+\frac{1}{2} (p + \frac{2}{n})-\frac{n-1}{n}\pi(\tau)-\frac{\alpha(\tau)}{n}$, is constant and the soliton is expanding, steady, shrinking according as $r-\phi-2(n-1)a+\frac{1}{2} (p + \frac{2}{n})-\frac{n-1}{n}\pi(\tau)-\frac{\alpha(\tau)}{n} \gtreqqless 0$. }\\\\
Now in \eqref{3.11}, if the 1-form $\alpha$ vanishes identically then $\lambda=r-\phi-2(n-1)a+\frac{1}{2} (p + \frac{2}{n})-\frac{n-1}{n}\pi(\tau).$\\\\
If the 1-form $\alpha$ vanishes identically and the function $\phi = 1$ in \eqref{3.11}, then $\lambda=r-1-2(n-1)a+\frac{1}{2} (p + \frac{2}{n})-\frac{n-1}{n}\pi(\tau).$\\\\
In \eqref{3.11}, if the function $\phi = 0$, then $\lambda=r-2(n-1)a+\frac{1}{2} (p + \frac{2}{n})-\frac{n-1}{n}\pi(\tau)-\frac{\alpha(\tau)}{n}.$\\\\
If $\phi = \alpha = 0$ in \eqref{3.11}, then $\lambda=r-2(n-1)a+\frac{1}{2} (p + \frac{2}{n})-\frac{n-1}{n}\pi(\tau).$\\\\
Finally in \eqref{3.11}, if $\alpha(\tau) = 0$, then $\lambda=r-\phi-2(n-1)a+\frac{1}{2} (p + \frac{2}{n})-\frac{n-1}{n}\pi(\tau).$\\
Then we have,\\\\
\textbf{Corollary 3.4.} {\em Let $(g,\tau,\lambda)$ be a conformal Yamabe soliton on $M$ with respect to the semi-symmetric metric connection $\bar{\nabla}$. Then the vector field $\tau$ is\\\\
(i)concircular if $\lambda=r-\phi-2(n-1)a+\frac{1}{2} (p + \frac{2}{n})-\frac{n-1}{n}\pi(\tau)$, is constant and the soliton is expanding, steady, shrinking according as $r-\phi-2(n-1)a+\frac{1}{2} (p + \frac{2}{n})-\frac{n-1}{n}\pi(\tau) \gtreqqless 0.$\\\\
(ii)concurrent if $\lambda=r-1-2(n-1)a+\frac{1}{2} (p + \frac{2}{n})-\frac{n-1}{n}\pi(\tau)$, is constant and the soliton is expanding, steady, shrinking according as $r-1-2(n-1)a+\frac{1}{2} (p + \frac{2}{n})-\frac{n-1}{n}\pi(\tau) \gtreqqless 0.$\\\\
(iii)recurrent if $\lambda=r-2(n-1)a+\frac{1}{2} (p + \frac{2}{n})-\frac{n-1}{n}\pi(\tau)-\frac{\alpha(\tau)}{n}$, is constant and the soliton is expanding, steady, shrinking according as $r-2(n-1)a+\frac{1}{2} (p + \frac{2}{n})-\frac{n-1}{n}\pi(\tau)-\frac{\alpha(\tau)}{n} \gtreqqless 0.$\\\\
(iv)parallel if $\lambda=r-2(n-1)a+\frac{1}{2} (p + \frac{2}{n})-\frac{n-1}{n}\pi(\tau)$, is constant and the soliton is expanding, steady, shrinking according as $r-2(n-1)a+\frac{1}{2} (p + \frac{2}{n})-\frac{n-1}{n}\pi(\tau) \gtreqqless 0.$\\\\
(v)torqued if $\lambda=r-\phi-2(n-1)a+\frac{1}{2} (p + \frac{2}{n})-\frac{n-1}{n}\pi(\tau)$, is constant and the soliton is expanding, steady, shrinking according as $r-\phi-2(n-1)a+\frac{1}{2} (p + \frac{2}{n})-\frac{n-1}{n}\pi(\tau) \gtreqqless 0.$}\\\\
Now we consider $(g,\tau,\lambda)$ as a conformal Yamabe soliton on $M$ with respect to the projective semi-symmetric connection $\tilde{\nabla}$.
Then we have,
\begin{equation}\label{3.12}
   (\tilde{\pounds}_\tau g)(X,Y)+[2\lambda -2\tilde{r} - (p + \frac{2}{n})]g(X,Y)=0,
\end{equation}
where $\tilde{\pounds}_\tau$ is the Lie derivative along $\tau$ with respect to $\tilde{\nabla}$.\\
Now from \eqref{2.5}, we have,
\begin{eqnarray}\label{3.13}
 (\tilde{\pounds}_\tau g)(X,Y) &=& g(\tilde{\nabla}_X \tau,Y)+g(X,\tilde{\nabla}_Y \tau) \nonumber\\
                               &=& g(\nabla_X \tau + \psi(\tau)X + \psi(X)\tau + \mu(\tau)X - \mu(X)\tau,Y) \nonumber\\
                               &+&g(X,\nabla_Y \tau + \psi(\tau)Y + \psi(Y)\tau + \mu(\tau)Y - \mu(Y)\tau)\nonumber\\
                               &=& (\pounds_\tau g)(X,Y)+\frac{1}{n+1}[2n \pi(\tau)g(X,Y)-\pi(X)g(\tau,X) \nonumber\\
                               &-&\pi(Y)g(\tau,X)].
\end{eqnarray}
Using \eqref{3.2} in \eqref{3.13}, we get,
\begin{eqnarray}\label{3.14}
 (\tilde{\pounds}_\tau g)(X,Y)&=&  2\phi g(X,Y)+\alpha(X)g(\tau,Y)+\alpha(Y)g(\tau,X) \nonumber \\
                              &+& \frac{1}{n+1}[2n \pi(\tau)g(X,Y)-\pi(X)g(\tau,X)-\pi(Y)g(\tau,X)].\nonumber\\
\end{eqnarray}
Now from \eqref{2.7} and \eqref{3.14}, \eqref{3.12} becomes,
\begin{multline}\label{3.15}
  \Big[\phi+\frac{n}{n+1}\pi(\tau)-r-\Tr(\theta)+(n-1)\Tr(\omega)+\lambda-\frac{1}{2}\Big(p + \frac{2}{n}\Big)\Big]g(X,Y) \\
  +\frac{1}{2}\Big[\Big\{\alpha(X)-\frac{\pi(X)}{n+1}\Big\}g(\tau,Y)+\Big\{\alpha(Y)-\frac{\pi(Y)}{n+1}\Big\}g(\tau,X)\Big]=0. \\
\end{multline}
Taking contraction of \eqref{3.15} over X and Y, we have,
\begin{equation}\label{3.16}
  \Big[\phi-r-\Tr(\theta)+(n-1)\Tr(\omega)+\lambda-\frac{1}{2}\Big(p + \frac{2}{n}\Big)\Big]n+(n-1)\pi(\tau)+\alpha(\tau)=0,
\end{equation}
which leads to,
\begin{equation}\label{3.17}
  \lambda=r-\phi+\Tr(\theta)-(n-1)\Tr(\omega)+\frac{1}{2}\Big(p + \frac{2}{n}\Big)-\frac{n-1}{n}\pi(\tau)-\frac{\alpha(\tau)}{n}.
\end{equation}
So we can state the following theorem:\\\\
\textbf{Theorem 3.5.} {\em Let $(g,\tau,\lambda)$ be a conformal Yamabe soliton on $M$ with respect to the projective semi-symmetric connection $\tilde{\nabla}$. Then the vector field $\tau$ is torse-forming if $ \lambda=r-\phi+\Tr(\theta)-(n-1)\Tr(\omega)+\frac{1}{2}(p + \frac{2}{n})-\frac{n-1}{n}\pi(\tau)-\frac{\alpha(\tau)}{n}$, is constant and the soliton is expanding, steady, shrinking according as $r-\phi+\Tr(\theta)-(n-1)\Tr(\omega)+\frac{1}{2}(p + \frac{2}{n})-\frac{n-1}{n}\pi(\tau)-\frac{\alpha(\tau)}{n} \gtreqqless 0$. }\\\\
Now in \eqref{3.17}, if the 1-form $\alpha$ vanishes identically then $ \lambda=r-\phi+\Tr(\theta)-(n-1)\Tr(\omega)+\frac{1}{2}(p + \frac{2}{n})-\frac{n-1}{n}\pi(\tau)$.\\\\
If the 1-form $\alpha$ vanishes identically and the function $\phi = 1$ in \eqref{3.17}, then $\lambda=r-1+\Tr(\theta)-(n-1)\Tr(\omega)+\frac{1}{2}(p + \frac{2}{n})-\frac{n-1}{n}\pi(\tau)$.\\\\
In \eqref{3.17}, if the function $\phi = 0$, then $ \lambda=r+\Tr(\theta)-(n-1)\Tr(\omega)+\frac{1}{2}(p + \frac{2}{n})-\frac{n-1}{n}\pi(\tau)-\frac{\alpha(\tau)}{n}$.\\\\
If $\phi = \alpha = 0$ in \eqref{3.17}, then $ \lambda=r+\Tr(\theta)-(n-1)\Tr(\omega)+\frac{1}{2}(p + \frac{2}{n})-\frac{n-1}{n}\pi(\tau)$.\\\\
Finally in \eqref{3.17}, if $\alpha(\tau) = 0$, then $ \lambda=r-\phi+\Tr(\theta)-(n-1)\Tr(\omega)+\frac{1}{2}(p + \frac{2}{n})-\frac{n-1}{n}\pi(\tau)$.\\
Then we have,\\\\
\textbf{Corollary 3.6.} {\em Let $(g,\tau,\lambda)$ be a conformal Yamabe soliton on $M$ with respect to the projective semi-symmetric connection $\tilde{\nabla}$. Then the vector field $\tau$ is\\\\
(i)concircular if $\lambda=r-\phi+\Tr(\theta)-(n-1)\Tr(\omega)+\frac{1}{2}(p + \frac{2}{n})-\frac{n-1}{n}\pi(\tau)$, is constant and the soliton is expanding, steady, shrinking according as $r-\phi+\Tr(\theta)-(n-1)\Tr(\omega)+\frac{1}{2}(p + \frac{2}{n})-\frac{n-1}{n}\pi(\tau) \gtreqqless 0.$\\\\
(ii)concurrent if $\lambda=r-1+\Tr(\theta)-(n-1)\Tr(\omega)+\frac{1}{2}(p + \frac{2}{n})-\frac{n-1}{n}\pi(\tau)$, is constant and the soliton is expanding, steady, shrinking according as $r-1+\Tr(\theta)-(n-1)\Tr(\omega)+\frac{1}{2}(p + \frac{2}{n})-\frac{n-1}{n}\pi(\tau) \gtreqqless 0.$\\\\
(iii)recurrent if $\lambda=r+\Tr(\theta)-(n-1)\Tr(\omega)+\frac{1}{2}(p + \frac{2}{n})-\frac{n-1}{n}\pi(\tau)-\frac{\alpha(\tau)}{n}$, is constant and the soliton is expanding, steady, shrinking according as $r+\Tr(\theta)-(n-1)\Tr(\omega)+\frac{1}{2}(p + \frac{2}{n})-\frac{n-1}{n}\pi(\tau)-\frac{\alpha(\tau)}{n} \gtreqqless 0.$\\\\
(iv)parallel if $\lambda=r+\Tr(\theta)-(n-1)\Tr(\omega)+\frac{1}{2}(p + \frac{2}{n})-\frac{n-1}{n}\pi(\tau)$, is constant and the soliton is expanding, steady, shrinking according as $r+\Tr(\theta)-(n-1)\Tr(\omega)+\frac{1}{2}(p + \frac{2}{n})-\frac{n-1}{n}\pi(\tau) \gtreqqless 0.$\\\\
(v)torqued if $\lambda=r-\phi+\Tr(\theta)-(n-1)\Tr(\omega)+\frac{1}{2}(p + \frac{2}{n})-\frac{n-1}{n}\pi(\tau)$, is constant and the soliton is expanding, steady, shrinking according as $r-\phi+\Tr(\theta)-(n-1)\Tr(\omega)+\frac{1}{2}(p + \frac{2}{n})-\frac{n-1}{n}\pi(\tau) \gtreqqless 0.$}\\
\vspace {0.3cm}
\section{\textbf{Application of torse forming vector field on $*$-Yamabe soliton}}
Let $(g,\tau,\lambda)$ be a $*$-Yamabe soliton on $M$ with respect to the Riemannian connection $\nabla$. Then from \eqref{1.7}, we get,
\begin{equation}\label{4.1}
  \frac{1}{2}(\pounds_\tau g)(X,Y) = (r^*-\lambda)g(X,Y).
\end{equation}
Using \eqref{3.2}, \eqref{4.1} becomes,
\begin{equation}\label{4.2}
  (r^*-\lambda-\phi)g(X,Y)=\frac{1}{2}[\alpha(X)g(\tau,Y)+\alpha(Y)g(\tau,X)].
\end{equation}
Taking contraction of \eqref{4.2} over X and Y, we have,
\begin{equation}\label{4.3}
  (r^*-\lambda-\phi)n=\alpha(\tau),
\end{equation}
which leads to,
\begin{equation}\label{4.4}
  \lambda=r^*-\phi-\frac{\alpha(\tau)}{n}.
\end{equation}
Hence we can state the following:\\\\
\textbf{Theorem 4.1} {\em Let $(g,\tau,\lambda)$ be a $*$-Yamabe soliton on $M$ with respect to the Riemannian connection $\nabla$. Then the vector field $\tau$ is torse-forming if $\lambda=r^*-\phi-\frac{\alpha(\tau)}{n}$, is constant and the soliton is expanding, steady, shrinking according as $r^*-\phi-\frac{\alpha(\tau)}{n} \gtreqqless 0$. }\\\\
Now in \eqref{4.4}, if the 1-form $\alpha$ vanishes identically then $\lambda=r^*-\phi.$\\
If the 1-form $\alpha$ vanishes identically and the function $\phi = 1$ in \eqref{4.4}, then $\lambda=r^*-1.$\\
In \eqref{4.4}, if the function $\phi = 0$, then $\lambda=r^*-\frac{\alpha(\tau)}{n}.$\\
If $\phi = \alpha = 0$ in \eqref{4.4}, then $\lambda=r^*.$\\
Finally in \eqref{4.4}, if $\alpha(\tau) = 0$, then $\lambda=r^*-\phi.$\\
Then we have,\\\\
\textbf{Corollary 4.2.} {\em Let $(g,\tau,\lambda)$ be a $*$-Yamabe soliton on $M$ with respect to the Riemannian connection $\nabla$. Then the vector field $\tau$ is\\\\
(i)concircular if $\lambda=r^*-\phi$, is constant and the soliton is expanding, steady, shrinking according as $r^*-\phi \gtreqqless 0$.\\\\
(ii)concurrent if $\lambda=r^*-1$, is constant and the soliton is expanding, steady, shrinking according as $r^*-1 \gtreqqless 0$.\\\\
(iii)recurrent if $\lambda=r^*-\frac{\alpha(\tau)}{n} $, is constant and the soliton is expanding, steady, shrinking according as $r^*-\frac{\alpha(\tau)}{n} \gtreqqless 0$.\\\\
(iv)parallel if $\lambda=r^*$, is constant and the soliton is expanding, steady, shrinking according as $r^* \gtreqqless 0$.\\\\
(v)torqued if $\lambda=r^*-\phi$, is constant and the soliton is expanding, steady, shrinking according as  $r^*-\phi \gtreqqless 0$.}\\
\vspace {0.3cm}
\section{\textbf{Example}}
Let $M = \{(x, y, z) \in \mathbb{R}^3, z \neq 0\}$ be a manifold of dimension 3, where $(x, y, z)$ are standard coordinates in $\mathbb{R}^3$. The vector fields,
\begin{equation}
  e_1 = z^2 \frac{\partial}{\partial x}, \quad e_2= z^2 \frac{\partial}{\partial y}, \quad e_3= \frac{\partial}{\partial z} \nonumber
\end{equation}
are linearly independent at each point of $M$.\\
Let $g$ be the Riemannian metric defined by
\begin{equation}
  g(e_1,e_2)=g(e_2,e_3)=g(e_3,e_1)=0,\nonumber
\end{equation}
\begin{equation}
   g(e_1,e_1) = g(e_2,e_2) = g(e_3,e_3) =1.\nonumber
\end{equation}
Let $\nabla$ be the Levi-Civita connection with respect to the Riemannian metric $g$. Then we have,
  $$ [e_1,e_2] =0, \quad [e_1,e_3] =-\frac{2}{z} e_1, \quad [e_2,e_3]= -\frac{2}{z} e_2.$$
The connection $\nabla$ of the metric $g$ is given by,
\begin{eqnarray}
  2g(\nabla_X Y,Z) &=& Xg(Y,Z)+Yg(Z,X)-Zg(X,Y)\nonumber \\
                   &-& g(X, [Y,Z])-g(Y, [X, Z]) + g(Z, [X, Y ]),\nonumber
\end{eqnarray}
which is known as Koszul’s formula.\\
Using Koszul’s formula, we can easily calculate,
$$\nabla_{e_1} e_1 =\frac{2}{z} e_3, \quad \nabla_{e_1} e_2 =0 ,\quad \nabla_{e_1} e_3 =- \frac{2}{z} e_1,$$
$$\nabla_{e_2} e_1 =0, \quad  \nabla_{e_2} e_2 =\frac{2}{z} e_3, \quad  \nabla_{e_2} e_3 =-\frac{2}{z} e_2,$$
$$\nabla_{e_3} e_1 =0, \quad  \nabla_{e_3} e_2 = 0, \quad  \nabla_{e_3} e_3 =0.$$
Also, the Riemannian curvature tensor $R$ is given by,
$$R(X, Y )Z = \nabla_X\nabla_Y Z - \nabla_Y \nabla_X Z - \nabla_{[X,Y]} Z.$$
Hence,
$$R(e_1,e_2)e_1 =\frac{4}{z^2} e_2, \quad R(e_1,e_2)e_2 = -\frac{4}{z^2} e_1, \quad R(e_1,e_3)e_1 =\frac{6}{z^2} e_3,$$
$$R(e_1,e_3)e_3 = -\frac{6}{z^2}e_1, \quad R(e_2,e_3)e_2 = \frac{6}{z^2}e_3, \quad R(e_2,e_3)e_3 = -\frac{6}{z^2}e_2.$$
$$R(e_1,e_2)e_3 = 0, \quad R(e_2,e_3)e_1 = 0, \quad R(e_3,e_1)e_2 = 0.$$
Then, the Ricci tensor $S$ is given by,
$$S(e_1,e_1) = -\frac{10}{z^2}, \quad S(e_2,e_2) = -\frac{10}{z^2}, \quad S(e_3,e_3)= -\frac{12}{z^2}.$$
Hence the scalar curvature is, $r=-\frac{32}{z^2}$.\\
Since $\{e_1, e_2, e_3\}$ forms a basis then any vector field $X, Y,W \in \chi(M)$ can be written as:
$$X = a_1e_1 + b_1e_2 + c_1e_3, \quad Y = a_2e_1 + b_2e_2 + c_2e_3,\quad W = a_3e_1 + b_3e_2 + c_3e_3,$$
where $a_i, b_i, c_i \in \mathbb{R}^+$ for $i = 1, 2, 3$ such that
$$\frac{a_1a_2 + b_1b_2}{c_1} + c_1\left(\frac{b_2}{b_1}-\frac{a_2}{a_1}- 1\right) \neq 0.$$
Now we choose the 1-form $\alpha$ by $\alpha(U)=g(U,\frac{2}{z}e_3)$ for any $U \in \chi(M)$ and the smooth function $\phi$ as:
$$\phi=\frac{2}{z}\left\{\frac{a_1a_2 + b_1b_2}{c_1} + c_1\left(\frac{b_2}{b_1}-\frac{a_2}{a_1}- 1\right)\right\}.$$
Then the relation
\begin{equation}\label{5.1}
  \nabla_X Y =\phi X +\alpha(X)Y
\end{equation}
holds. Hence $Y$ is a torse-forming vector field.\\
Then from \eqref{5.1}, we obtain,
\begin{eqnarray}\label{5.2}
  (\pounds_Y g)(X, W) &=&  g(\nabla_X Y,W)+g(X,\nabla_W Y) \nonumber\\
                        &=& 2\phi g(X,W)+\alpha(X)g(Y,W)+\alpha(W)g(Y,X).
\end{eqnarray}
Also we have,
\begin{eqnarray}\label{5.3}
 g(X,Y)&=&a_1a_2 + b_1b_2 + c_1c_2, \nonumber\\
  g(Y,W)&=&a_2a_3 + b_2b_3 + c_2c_3, \nonumber\\
  g(X,W)&=&a_1a_3 + b_1b_3 + c_1c_3,
\end{eqnarray}
and
\begin{equation}\label{5.4}
  \alpha(X)=\frac{2c_1}{z}, \quad \alpha(Y)= \frac{2c_2}{z}, \quad \alpha(W)=\frac{2c_3}{z}.
\end{equation}
From \eqref{5.3} and \eqref{5.4}, \eqref{5.2} becomes,
\begin{eqnarray}\label{5.5}
  (\pounds_Y g)(X, W) &=& \frac{2}{z}\Big[\Big\{\frac{2(a_1a_2 + b_1b_2)}{c_1} + 2c_1\Big(\frac{b_2}{b_1}-\frac{a_2}{a_1}- 1\Big)\Big\}(a_1a_3 + b_1b_3 + c_1c_3) \nonumber\\
                      &+& c_1(a_2a_3 + b_2b_3 + c_2c_3)+c_3(a_1a_2 + b_1b_2 + c_1c_2)\Big],
\end{eqnarray}
and
\begin{equation}\label{5.6}
  \Big[2\lambda-2r-\Big(p + \frac{2}{3}\Big)\Big]g(X,W)=2\Big[\lambda+\frac{32}{z^2}-\frac{1}{2}\Big(p + \frac{2}{3}\Big)\Big](a_1a_3 + b_1b_3 + c_1c_3).
\end{equation}
Let us assume that $a_1a_3 + b_1b_3 + c_1c_3 \neq 0$ and
\begin{equation}\label{5.7}
  3c_1(a_2a_3 + b_2b_3 + c_2c_3)+3c_3(a_1a_2 + b_1b_2 + c_1c_2)-2c_2(a_1a_3 + b_1b_3 + c_1c_3)=0.
\end{equation}
Hence $(g, Y, \lambda)$ is a conformal Yamabe soliton on $M$, i.e.
$$(\pounds_Y g)(X, W)-2rg(X,W)+\Big[2\lambda-\Big(p + \frac{2}{3}\Big)\Big]g(X,W) = 0,$$
provided,
\begin{eqnarray}
  \lambda &=& -\frac{32}{z^2}-\frac{2}{z}\Big\{\frac{a_1a_2 + b_1b_2}{c_1} + c_1\Big(\frac{b_2}{b_1}-\frac{a_2}{a_1}- 1\Big)\Big\} \nonumber\\
          &-& \frac{c_1(a_2a_3 + b_2b_3 + c_2c_3)+c_3(a_1a_2 + b_1b_2 + c_1c_2)}{z(a_1a_3 + b_1b_3 + c_1c_3)}+\frac{1}{2}\Big(p + \frac{2}{3}\Big)\nonumber \\
          &=& r-\phi-\frac{1}{3}\alpha(Y)+\frac{1}{2}\Big(p + \frac{2}{3}\Big)\quad (using ~ \eqref{5.4} ~ and ~ \eqref{5.7}) \nonumber \\
          &=& constant.\nonumber
\end{eqnarray}
Hence the condition of existence of the conformal Yamabe soliton $(g,Y,\lambda)$ on a 3-dimensional Riemannian manifold $M$ with potential vector field $Y$ as torse forming in Theorem:3.1 is satisfied.\\

\end{document}